\documentclass[10pt]{amsart}

\usepackage[nobysame,abbrev,alphabetic]{amsrefs}
\usepackage{amssymb}
\usepackage[all]{xy}
\usepackage{verbatim}

\DeclareMathOperator{\cd}{cd}
\DeclareMathOperator{\Char}{char}
\DeclareMathOperator{\Gal}{Gal}

\DeclareMathOperator{\Sym}{Sym}

\begin{document}

\newtheorem{thm}{Theorem}[section]
\newtheorem{cor}[thm]{Corollary}
\newtheorem{lem}[thm]{Lemma}
\newtheorem{prop}[thm]{Proposition}
\newtheorem{defin}[thm]{Definition}
\newtheorem{exam}[thm]{Example}
\newtheorem{examples}[thm]{Examples}
\newtheorem{rem}[thm]{Remark}
\newtheorem{case}{\sl Case}
\newtheorem{claim}{Claim}
\newtheorem{prt}{Part}
\newtheorem*{mainthm}{Main Theorem}
\newtheorem*{thmA}{Theorem A}
\newtheorem*{thmB}{Theorem B}
\newtheorem*{thmC}{Theorem C}
\newtheorem*{thmD}{Theorem D}
\newtheorem*{thmE}{Theorem E}
\newtheorem{question}[thm]{Question}
\newtheorem*{notation}{Notation}
\swapnumbers
\newtheorem{rems}[thm]{Remarks}
\newtheorem*{acknowledgment}{Acknowledgment}

\newtheorem{questions}[thm]{Questions}
\numberwithin{equation}{section}

\newcommand{\ab}{\mathrm{ab}}
\newcommand{\Coker}{\mathrm{Coker}}
\newcommand{\dec}{\mathrm{dec}}
\newcommand{\dirlim}{\varinjlim}
\newcommand{\discup}{\ \ensuremath{\mathaccent\cdot\cup}}
\newcommand{\gr}{\mathfrak{gr}}
\newcommand{\nek}{,\ldots,}
\newcommand{\inv}{^{-1}}
\newcommand{\isom}{\cong}
\newcommand{\sep}{\mathrm{sep}}
\newcommand{\sym}{\mathrm{sym}}
\newcommand{\tagg}{^{''}}
\newcommand{\tensor}{\otimes}

\newcommand{\alp}{\alpha}
\newcommand{\gam}{\gamma}
\newcommand{\Gam}{\Gamma}
\newcommand{\del}{\delta}
\newcommand{\Del}{\Delta}
\newcommand{\eps}{\epsilon}
\newcommand{\lam}{\lambda}
\newcommand{\Lam}{\Lambda}
\newcommand{\sig}{\sigma}
\newcommand{\Sig}{\Sigma}
\newcommand{\bfA}{\mathbf{A}}
\newcommand{\bfB}{\mathbf{B}}
\newcommand{\bfC}{\mathbf{C}}
\newcommand{\bfF}{\mathbf{F}}
\newcommand{\bfP}{\mathbf{P}}
\newcommand{\bfQ}{\mathbf{Q}}
\newcommand{\bfR}{\mathbf{R}}
\newcommand{\bfS}{\mathbfS}
\newcommand{\bfT}{\mathbf{T}}
\newcommand{\bfZ}{\mathbf{Z}}
\newcommand{\dbA}{\mathbb{A}}
\newcommand{\dbC}{\mathbb{C}}
\newcommand{\dbF}{\mathbb{F}}
\newcommand{\dbN}{\mathbb{N}}
\newcommand{\dbQ}{\mathbb{Q}}
\newcommand{\dbR}{\mathbb{R}}
\newcommand{\dbZ}{\mathbb{Z}}
\newcommand{\grf}{\mathfrak{f}}
\newcommand{\gra}{\mathfrak{a}}
\newcommand{\grm}{\mathfrak{m}}
\newcommand{\grp}{\mathfrak{p}}
\newcommand{\grq}{\mathfrak{q}}
\newcommand{\calA}{\mathcal{A}}
\newcommand{\calB}{\mathcal{B}}
\newcommand{\calC}{\mathcal{C}}
\newcommand{\calE}{\mathcal{E}}
\newcommand{\calG}{\mathcal{G}}
\newcommand{\calH}{\mathcal{H}}
\newcommand{\calK}{\mathcal{K}}
\newcommand{\calL}{\mathcal{L}}
\newcommand{\calW}{\mathcal{W}}
\newcommand{\calV}{\mathcal{V}}

\title[A small quotient of the big absolute Galois group]{A small quotient of the big absolute Galois group (joint work with Ido Efrat and J\'{a}n Min\'{a}\v{c})}

\author{Sunil K.\ Chebolu}
\address{Department of Mathematics\\
Illinois State University\\
Campus box 4520\\
Normal, IL 61790\\
USA} \email{schebol@ilstu.edu}

\maketitle

\section{Introduction}
This is a report of a talk given at the Oberwolfach workshop on ``cohomology of finite groups: Interactions and applications''
which was held  during July 25th - July 31st, 2010.
It is an announcement of some of the results (with motivation) and their applications from the paper ``Quotients of absolute Galois groups which determine the entire Galois cohomology''; see arXiv:0905.1364.

Let $F$ be a field and let $F_{\sep}$ denote the separable closure of $F$. The absolute Galois group $G_F$ of $F$ is the Galois group of the Galois 
extention $F_{\sep}/F$. It can be expressed as the inverse limit of all the finite Galois extentions over $F$:
\[ G_F:= \text{Gal}(F_{\sep}/F) =  \underset{[L \colon F] < \infty}{\text{lim}} \text{Gal} (L/F).  \]
Thus, absolute Galois groups are profinite groups. They are very rich but mysterious objects which are of great interest. For instance,
they play a central role in study of the inverse Galois problem: given a finite group $G$ when it is realisable as a Galois group of some Galois extension over $L$ over $F$? note that from the Galois correspondence, this question has an affirmative answer provided there is a normal quotient of the absolute Galois group $G_F$ by some
closed normal subgroup whose quotient is isomorphic to $G$. Therefore a good knowledge of the structure of the absolute Galois groups is very
desirable. The inverse Galois problem is wide open. It is known in some cases. For example, the Kronecker-Weber theorem tells us that finite abelian groups are realisable over $\mathbb{Q}$, the field of rational numbers. Similary, a deeper result of Shafarevich states that every solvable finite group is realisable over $\mathbb{Q}$.

 A main open problem in modern Galois theory is the characterization of the
profinite groups which are realizable as absolute Galois groups of fields $F$.
The non-trivial torsion in such groups is described by the Artin--Schreier theory from the late 1920's,
namely, it consists solely of involutions. That is, the only non-trivial finite subgroup of an absolute Galois group is $C_2$, the
cyclic group of order $2$.
More refined information on the structure of absolute Galois groups is given by Galois cohomology,
systematically developed starting the 1950's by Tate, Serre, and others.
Yet, explicit examples of torsion-free profinite groups which are not absolute Galois groups are rare.
In 1970, Milnor  \cite{Mil70} introduced his $K$-ring functor $K^M_*(F)$, and pointed out
close connections between this graded ring and the mod-$2$ Galois cohomology of the field.
This connection, in a more general form, became  known as the Bloch--Kato conjecture:
it says that for all $r\geq0$ and all $m$ prime to  $\Char\, F$,
there is a canonical isomorphism $K^M_r(F)/m\to H^r(G_F,\mu_m^{\tensor r})$.
The conjecture was proved for $r=2$ by Merkurjev and Suslin \cite{MerkurjevSuslin82},
for $r$ arbitrary and $m=2$ by Voevodsky \cite{Voevodsky03a},
and in general by Rost, Voevodsky, with a patch by Weibel (\cite{Voevodsky03b}, \cite{Weibel09},
\cite{Weibel08}). In particular, the Block-Kato conjecture implies that the Galois cohomology ring has generators in degree one and
relations in degree two. 

We obtain new constrains on the group structure of
absolute Galois groups of fields, using this isomorphism.
We use these constrains to produce new examples of torsion-free profinite groups which are not absolute Galois groups.
We also demonstrate that the maximal pro-$p$ quotient of the absolute
Galois group can be characterized in purely cohomological terms!
The main object of the talk is a remarkable small quotient of the absolute Galois group,
which, because of the above isomorphism, already carries a substantial information about the arithmetic of $F$.

\section{Main theorems}
In this section we mention a couple of results. For more results and proofs of the results we refer the reader to our paper 
on the arXiv: 0905.1364. 
Fix a prime number $p$ and a $p$-power $q=p^d$, with $d\geq1$.
All fields which appear in this paper will be tacitly assumed to contain a primitive $q$th root of unity.
Let $F$ be such a field and let $G_F=\Gal(F_\sep/F)$ be its absolute Galois group, where $F_\sep$ is
the separable closure of $F$.
Let $H^*(G_F)=H^*(G_F,\dbZ/q)$ be the Galois cohomology ring with the trivial action of $G_F$ on $\dbZ/q$.
Our new constraints relate the descending $q$-central sequence $G_F^{(i)}$,
$i=1,2,3\nek$ of $G_F$  with $H^*(G_F)$.
Setting $G_F^{[i]}=G_F/G_F^{(i)}$, we show that the quotient $G_F^{[3]}$ determines $H^*(G_F)$, and vice versa.
Specifically, we prove:

\begin{thmA}
The inflation map gives an isomorphism
\[
H^*(G_F^{[3]})_\dec\xrightarrow{\sim} H^*(G_F),
\]
where $H^*(G_F^{[3]})_\dec$ is the decomposable part of $H^*(G_F^{[3]})$
(i.e., its subring generated by degree $1$ elements).
\end{thmA}

We further establish the following result.

\begin{thmB}
Let $F_1$, $F_2$ be fields and let $\pi\colon G_{F_1}\to G_{F_2}$ be a
(continuous) homomorphism.
The following conditions are equivalent:
\begin{enumerate}
\item[(i)]
the induced map $\pi^*\colon H^*(G_{F_2})\to H^*(G_{F_1})$ is an isomorphism;
\item[(ii)]
the induced map $\pi^{[3]}\colon G_{F_1}^{[3]}\to G_{F_2}^{[3]}$ is an
isomorphism.
\end{enumerate}
\end{thmB}

Theorems A and B show that $G_F^{[3]}$ is a Galois-theoretic analog of the
cohomology ring $H^*(G_F)$.
Its structure is considerably simpler and more accessible than the full absolute Galois group $G_F$
(see e.g., \cite{EfratMinac}).
Yet, as shown in our theorems, these small and accessible quotients encode and control the entire
cohomology ring.

In the case $q=2$ the group $G_F^{[3]}$ has been extensively studied under the name ``$W$-group",
in particular in connection with quadratic forms \cite{MinacSpira96},
\cite{AdemKaraMinac99}, \cite{MaheMinacSmith04}). In this special case,
Theorem A was proved in \cite{AdemKaraMinac99}*{Th.\ 3.14}. It was further
shown that then $G_F^{[3]}$ has great arithmetical significance: it encodes
large parts of the arithmetical structure of $F$, such as its orderings, its
Witt ring, and certain non-trivial valuations. Theorem A explains this
surprising phenomena, as these arithmetical objects are known to be encoded in
$H^*(G_F)$ (with the additional knowledge of the Kummer element of $-1$).

First links between these quotients and the Bloch--Kato conjecture, and its special case the Merkurjev--Suslin theorem,
were already noticed in a joint work of Min\'a\v c and Spira and in work of  Bogomolov.

Our approach is purely group-theoretic, and the main results above are in fact
proved for arbitrary profinite groups which satisfy certain conditions on
their cohomology.
A key point is a rather general group-theoretic approach, partly inspired by \cite{GaoMinac97},
to the Milnor $K$-ring construction by means of quadratic hulls
of graded algebras.
The Rost--Voevodsky theorem on the bijectivity of the Galois symbol shows that these
cohomological conditions are satisfied by absolute Galois groups as above.
Using this we deduce aforementioned Theorems  in their field-theoretic version.

\section{Groups which are not maximal pro-$p$ Galois groups}
\label{section on examples}
We now apply Theorem A to give examples of pro-$p$ groups which cannot be realized
as maximal pro-$p$ Galois groups of fields (assumed as before to contain a root of unity of order $p$).
The groups we construct are only a sample of the most simple and straightforward examples illustrating our theorems, and
many other more complicated examples can be constructed along the same lines.
Throughout this section $q=p$. We have the following immediate consequence of
the analog of Theorem A for maximal pro-$p$ Galois groups:

\begin{prop}
\label{principle}
If $G_1,G_2$ are pro-$p$ groups such that $G_1^{[3,p]}\isom G_2^{[3,p]}$ and that
$H^*(G_1)\not\isom H^*(G_2)$, then at most one of them can be realized as the maximal pro-$p$ Galois group of a field.
\end{prop}

\begin{cor}
\label{cor to principle}
Let $S$ be a free pro-$p$ group and $R$ a nontrivial closed normal subgroup of $S^{(3,p)}$.
Then $G=S/R$ cannot occur as a maximal pro-$p$ Galois group of a field.
\end{cor}

\begin{exam}  \rm
Let $S$ be a free pro-$p$ group on $2$ generators, and take $R=[S,[S,S]]$.
By Corollary \ref{cor to principle}, $G=S/R$ is not realizable as $G_F(p)$ for a field $F$ as above.
Note that $G/[G,G]\isom S/[S,S]\isom\dbZ_p^2$ and $[G,G]=[S,S]/S^{(3,0)}\isom\dbZ_p$, so $G$ is torsion-free.
\end{exam}

\begin{prop}
\label{cohomological dimension}
Let $G$ be a pro-$p$ group such that $\dim_{\dbF_p}H^1(G)<\cd(G)$.
When $p=2$ assume also that $G$ is torsion-free.
Then $G$ is not a maximal pro-$p$ Galois group of a field as above.
\end{prop}

\begin{exam}  \rm
\label{wreath}
Let $K,L$ be finitely generated pro-$p$ groups with $1\leq n=\cd(K)<\infty$, $\cd(L)<\infty$, and $H^n(K)$ finite.
Let $\pi\colon L\to \Sym_m$, $x\mapsto \pi_x$,
be a homomorphism such that $\pi(L)$ is a transitive subgroup of $\Sym_m$.
Then $L$ acts on $K^m$ from the left by ${}^x(y_1\nek y_m)=(y_{\pi_x(1)}\nek y_{\pi_x(m)})$.
Let $G=K^m\rtimes L$.
It is generated by the generators of one copy of $K$ and of $L$.
Hence $\dim_{\dbF_p}H^1(G)=\dim_{\dbF_p}H^1(K)+\dim_{\dbF_p}H^1(L)$.

On the other hand, a routine inductive spectral sequence argument (see \cite{NeukirchSchmidtWingberg}*{Prop.\ 3.3.8})
shows that for every $i\geq0$ one has
\begin{enumerate}
\item
$\cd(K^i)=in$;
\item
$H^{in}(K^i)=H^n(K,H^{(i-1)n}(K^{i-1}))$, with the trivial $K$-action, is finite.
\end{enumerate}
Moreover, $\cd(G)=\cd(K^m)+\cd(L)$.
For $m$ sufficiently large we get $\dim_{\dbF_p}H^1(G)<mn+\cd(L)=\cd(G)$,
so by Proposition \ref{cohomological dimension}, $G$ is not a maximal pro-$p$ Galois group as above.
When $K,L$ are torsion-free, so is $G$.

For instance, one can take $K$ to be a free pro-$p$ group $\neq1$ on finitely many generators,
and let $L=\dbZ_p$ act on the direct product of $p^s$ copies of $K$ via $\dbZ_p\to\dbZ/p^s$
by cyclicly permuting the coordinates.
\end{exam}

\end{document}